\documentstyle[11pt]{article}

\newcommand{\comment}[1]{}

\newtheorem{th}{Theorem}[section]

\newtheorem{Theorem}[th]{Theorem}
\newtheorem{Lemma}[th]{Lemma}
\newtheorem{Proposition}[th]{Proposition}
\newtheorem{Remark}[th]{Remark}

\newtheorem{ccote}[th]{}

\newcommand{\decale}[1]{\par\noindent\hskip 3em\llap{#1\enspace}\ignorespaces}
\newcommand{\SS}{\S \kern .2em}

\newcommand{\preu}{{\sc Proof: \ }}

\newcommand{\hfl}[2]{\smash{\mathop{\hbox to 1 truecm{\kern 3pt\rightarrowfill\kern 3pt}}%
\limits^{\scriptstyle#1}_{\scriptstyle#2}}}
\newcommand{\cqfd}{\unskip\kern 6pt\penalty 500
\raise -2pt\hbox{\vrule\vbox to10pt{\hrule width
 4pt\vfill\hrule}\vrule}\smallskip}
\newcommand{\eqref}[1]{(\ref{#1})}
\newcommand{\bbr}{{\bf R}}
\newcommand{\bbp}{{\bf P}}
\newcommand{\bbc}{{\bf C}}
\newcommand{\bbz}{{\bf Z}}

\newcommand{\cala}{{\cal A}}

\newcommand{\cali}{{\cal I}}
\newcommand{\calj}{{\cal J}}
\newcommand{\calk}{{\cal K}}
\newcommand{\call}{{\cal L}}
\newcommand{\calm}{{\cal M}}
\newcommand{\caln}{{\cal N}}
\newcommand{\calo}{{\cal O}}

\newcommand{\cals}{{\cal S}}

\newcommand{\pcirc}{\kern .7pt {\scriptstyle \circ} \kern 1pt}

\newcommand{\pol}{{\rm Pol\,}}
\newcommand{\good}{lopsided}

\newcommand{\modmod}[1]{\mathop{\,\big/\!\!\big/\,}_{\!\!\!\textstyle #1}}

\title{Maximal Hamiltonian tori for polygon spaces}
\author{Jean-Claude HAUSMANN \and Susan TOLMAN
\footnote{Both authors thank the Swiss National Fund 
for Scientific  Research for its support.
The second author is partially supported
by a Sloan Fellowship and a National Science Foundation Grant.}}
\date{June 21 2002}

\begin{document} \maketitle

\begin{abstract}
We study the poset of Hamiltonian tori for polygon spaces.
We determine some maximal elements and give examples
where maximal Hamiltonian tori are not all of the same 
dimension.
\end{abstract}

\section{Introduction}\label{intro}

Let $M$ be a symplectic manifold and let $\cals (M)$ be the group
of symplectomorphisms of $M$. A sub-torus of $\cals (M)$ is called a 
{\it symplectic torus}; these tori are partially ordered by inclusions. 
In this paper, we study the maximal symplectic tori of  
polygon spaces with a particular emphasis on bending tori
(see the definitions below). Since polygon spaces are simply connected,
symplectic tori act on $M$ in a Hamiltonian fashion 
so we refer to them as {\it Hamiltonian tori}.

Let $E$ be a finite set together with a function $\lambda : E \to \bbr_+$. 
Define the space $\widetilde \pol (E,\lambda)$ by
$$\widetilde\pol(E,\lambda):=\bigg\{\rho : E\to\bbr^3  \biggm |
\sum_{e\in E}\rho (e)=0 \hbox{ and }  |\rho(e)|=\lambda(e)\ \forall e\in E
\bigg\} \ .$$
The {\it polygon space} $\pol (E,\lambda)$ is the quotient
$\pol(E,\lambda):=\widetilde\pol(E,\lambda)\bigm/SO_3$.
By choosing a bijection between $E$ and $\{1,\dots ,m\}$, the space
$\pol (E,\lambda)$ is regarded as the space of configurations in
$\bbr^3$ of a polygon with $m$ edges of length
$\lambda _1,\dots ,\lambda _m$, modulo rotation, whence the name
``polygon space". Also, we call
an element of $E$ an {\it edge} and $\lambda$ the {\it length function}.

A length function $\lambda $ is called {\it generic} if there is no map
$\varepsilon : E\to \{\pm 1\}$ so that
$\sum_{e\in E}\varepsilon (e)\lambda (e)=0$. 
This guarantees that the polygon cannot collapse to a line. 
In this paper, we always assume that $\lambda $ is generic and 
that $\pol (E,\lambda )$ is not empty.
In this case, $\pol (E,\lambda )$ is a closed smooth symplectic
manifold of dimension $2(|E|-3)\geq 0$. The polygon spaces are
better known as the moduli spaces of (weighted) ordered points on $\bbp^ 1$,
and also arise via other symplectic reductions (see [Kl], [KM], [HK1] and 
the proof of Proposition \ref{redpoly} below).

A subset $I$ of $E$ is called {\it \good\ } if there exists $e_0\in I$ such that
$\lambda (e_0) > \sum_{e\in I - \{e_0\}}\lambda (e)$.
The empty set is not lopsided, while
a singleton $\{ e\}$ is always \good\ since the length function
takes strictly positive values. The total set $E$ is not \good\ since
$\pol(E,\lambda )$ is assumed to be non-empty.

For $I\subset E$ define $\rho_I : \widetilde\pol(E,\lambda )\to \bbr^3$ 
by $\rho_I:=\sum_{e\in I}\rho(e)$. The continuous function 
and $f_I : \widetilde\pol(E,\lambda )\to \bbr$ by
$f_I(\rho):=|\sum_{i\in I}\rho_i|$ descends to a function on
$\pol(E,\lambda )$, still called $f_I$. When $I$ is 
\good , this function 
does not vanish and is therefore smooth. Its Hamiltonian flow
$\Phi_I^t$ is called the {\it bending flow} associated to $I$. 
Bending flows have been introduced in [Kl] and [KM]. 
They are periodic (see [Kl, \SS 2.1] or [KM, Corollary 3.9]): 
$\Phi_I^t$ rotates at constant speed 
the set of vectors $\{\rho(e) \mid e\in I\}$ around the axis $\rho_I$.

A {\it bending torus} is a  Hamiltonian torus in
$\cals(\pol(E,\lambda ))$ generated by bending flows. Since the dimension of
$\pol(E,\lambda )$ is $2(|E|-3)$, the dimension of any Hamiltonian torus
is  at most  $|E|-3$.

In this paper, we study the poset of bending tori and compare it
with that of Hamiltonian ones.
For instance, the  following result is proved in Section \ref{mbt}
(see Corollary  \ref{th1}):  

\smallskip \noindent
{\bf Theorem A}\ \sl  Let $N(\lambda )$ be the minimal number of
\good\ subsets which are  necessary for a partition of $E$.
Then the maximal dimension of a bending torus for $\pol (E,\lambda)$ is 
$|E|-\max \{3,N(\lambda )\}$.\vskip .2 truecm \rm

We also give a more general statement that allows us to 
characterize maximal bending tori.
In some cases, these coincide with maximal Hamiltonian tori:

\smallskip \noindent
{\bf Theorem B}\ \sl  
Let $T$ be a bending torus of $\pol (E,\lambda)$ of
dimension $\geq |E|-5$. Then $T$ is a maximal Hamiltonian torus if and only if  it is a
maximal bending torus. \vskip .2 truecm \rm

In Section \ref{ex}, we give several examples
where maximal Hamiltonian tori are not all of
the same dimension.
Using the work of Y. Karshon \cite{Ka},
we show the existence of
Hamiltonian tori which are not conjugate to
a bending torus (Proposition \ref{pnotbend}).
Finally, the relationship with
maximal tori in the contactomorphism group
of pre-quantum circle bundles, due to E. Lerman
\cite{Le}, is mentioned in \ref{Lerman}.

\section{Preliminaries - Bending sets}  \label{prel}

\begin{Lemma} \label{commu} Let $\cali$ be a family of
\good\ subsets of $E$. The following
conditions are equivalent:
\decale{a)} The bending flows  $\{\Phi_{I}^t\mid I\in \cali\}$ generate a
bending torus.
\decale{b)} For each pair $A,B  \subseteq  \cali$,  
either $A\cap B=\emptyset$ or one is contained into
the other.
\end{Lemma}

\preu By [Kl, \SS 2.1] or [KM, Corollary 3.9], the bending flows are
periodic. Therefore, a) is equivalent to the fact
that $\{f_A,f_B\}=0$ for all $A,B\in\cali$, where $\{\cdot ,\cdot\}$ denotes the Poisson
bracket. Proposition 2.1.2 of [Kl] shows that $\{f^2_A,f^2_B\}=0$ if and only if  the pair $A,B$
satisfies Condition b). Since $f_A$ and $f_B$ never vanish, the formula
$$ \{f^2_A,f^2_B\} = 4\,f_Af_B\,\{f_A,f_B\}$$
implies that $\{f^2_A,f^2_B\}=0$ if and only if $\{f_A,f_B\}=0$.  \cqfd   \smallskip

A set $\cali$ of lopsided subsets of $E$ is called a {\it bending set} 
if it contains every singleton $\{e\}$ and satisfies 
the following ``absorption condition":
\sl for each pair $A,B  \subseteq  \cali$,  either $A\cap B=\emptyset$ 
or one is contained in the other.\rm
\smallskip

Bending sets are  technically convenient to parametrize bending tori.
Indeed, let $\cali$ be a bending set. By \ref{commu}, the bending flows
$\{\Phi_{I}^t\mid I\in \cali\}$ generate a bending torus $T_\cali$.
Conversely, if $T$ is a bending torus, there is at least one set $\cali$ of
\good\  subsets satisfying the absorption condition such that $T=T_\cali$, 
and one can add singletons to $\cali$ to make it a bending set. 

The elements of $\cali$ are
partially ordered by inclusions, so one can 
associate to $\cali$ the family 
$\calm_\cali$ of its maximal elements.
A direct consequence of the definition is 
that $\calm_\cali$ is a partition of $E$.

A bending set $\cali$ is called {\it full} if, for each  $I\in\cali$ 
which is not a singleton, 
there exist $I',I''\in\cali$ so that $I$ is the disjoint union of $I'$ and $I''$. 
It is easy to check that 
this condition is equivalent to either of the following.
\decale{a)} Given $I$ and  $I'$ in $\cali$ such that $I'\subset I$, the union
$\cali \cup \{I'\}$ is not a bending set. This  justifies
 the term ``full": one can no longer add elements to 
$\cali$ and keep the latter a bending set.
\decale{b)} For all $I\in\cali$ the set 
$\{I'\in \cali : I'\subseteq I\}$ contains $2\,| I| -1$ elements.

\paragraph{Remark} Let $\cali$ be a bending set. The reader might 
find it helpful to consider
the graph of this poset. It is a union of disjoint trees, each of which 
contains a unique maximal element. The bending set $\cali$ is full 
iff these trees are binary: each vertex has one edge leaving it (except
the maximal ones which have none) and 2 edges pointing into it (except the singletons
which have none).

\begin{Lemma} \label{filling} Let $\cali$ be a bending set. Then
there exists a (non-unique) bending set $\hat\cali$ such that the following
conditions hold
\decale{1)} $  \cali\subset \hat\cali$ (therefore $T_\cali \subset
T_{\hat\cali}$).
\decale{2)} $\hat\cali$ is full.
\decale{3)} $\calm_{\hat\cali}=\calm_{\cali}$.
\end{Lemma}

\preu  
If $\cali$ is full we are done. Otherwise, we proceed by induction on the 
number of ``non-full'' elements of $\cali$: those 
$I\in\cali$ which are
not singletons and are not the disjoint union of 2 elements of $\cali$. 
Let $I\in\cali$ be  a minimal ``non-full'' element.

Let $I_1,\dots ,I_r$ be the maximal proper subsets of $I$ which are elements
of $\cali$. One of them, say $I_1$, contains the longest edge of $I$.
For $i=2,\dots ,r-1$, define $R_i:= I_1\cup \cdots \cup I_i$ and let
$\check \cali :=\cali\cup \{R_2\}\cup \cdots \cup \{R_{r-1}\}$.  
One has $I=R_{r-1}\sqcup I_r$, $R_{r-1}=R_{r-2}\sqcup I_{r-1}$ etc. 
As $I$ was minimal, it is no longer non-full in 
$\hat \cali$. This gives the inductive step. \cqfd

We shall now compute the dimension of a bending tori. We need some knowledge 
about the critical points of the maps $f_I$ and its symplectic reduction.
The following lemma comes from \cite[Theorem 3.2]{Ha}.

\begin{Lemma}\label{critfi}
Let $I$ be a lopsided subset of $E$. An element $\rho\in \pol (E,\lambda)$ is a critical
point for $f_I$ if and only if either the set $\{\rho(e)\mid e\in I \}$ or the set
$\{\rho(e)\mid e\notin I \}$ lies in a line. \cqfd
\end{Lemma}

\begin{Proposition}\label{redpoly} 
Let $A\subset E$.  
Define $\bar A := A\cup\{A\}$ and $\lambda^{A,t} : \bar A\to\bbr$
by $\lambda^{A,t}(e):=\lambda (e)$ for $e\in A$ and $\lambda^{A,t}(A):=t$.
Then, if $A$ is lopsided, the symplectic reduction of $\pol (E,\lambda)$ at $t$, for the action
of the bending circle $T_A$, is symplectomorphic to the 
product of the two polygon spaces
$$\pol (E,\lambda)\modmod{t} T_A\ \cong\ \pol (\bar A,\lambda^{A,t})\times 
\pol (\overline{E-A},\lambda^{E-A,t}).$$
\end{Proposition}

\begin{Remark}\label{redpolyrem} \rm  
a) Proposition \ref{redpoly} holds true 
even if $t$ is not a regular value. If it is,
the two right hand polygon spaces 
of the formula are generic by Lemma \ref{critfi}.
\vskip .2 truecm
b) The following is clear from the proof below: if $T_\cali$ 
is a bending torus and $A\in \cali$, 
then the action of $T_\cali$ descends to the reduced
space, giving rise to a product of two bending tori: one for the bending set
$\{I\in\cali\mid I\subset A\}$ and the other for $\{I\in\cali\mid I\not\subset A\}$
\vskip .2 truecm
c) In this paper, Proposition \ref{redpoly} is used only for $|A|=2$. In this case,
the reduction of $\pol (E,\lambda)$ at $t$ is symplectomorphic 
to a polygon space with $|E|-1$ 
edges, since $\pol (\bar A,\lambda^{A,t})$ is a point. However, the 
hypothesis $|A|=2$ does not simplify the proof.
\end{Remark}

\noindent{\sc Proof of Proposition \ref{redpoly} : } First recall the precise definition
for the symplectic structure on $\pol(E,\lambda)$ (for details, see \cite[\S\ 1]{HK1}).
For $s\in\bbr$, let $\calo (s)$ the coadjoint orbit of $SO(3)$ with symplectic volume 
$2s$. With the usual identification of $so(3)^*$ with $\bbr^3$, 
$\calo (s)$ is the 2-sphere centered in $0$ of radius $r$.
For $A\subset E$, let $\mu_A : \prod_{e\in E}\calo(\lambda(e))\to\bbr^3$
be the partial sum $\mu_A((z_e)):=\sum_{e\in A}z_e$.
This is the moment map for the diagonal action of $SO(3)$ 
on the component indexed by $e\in A$. 
The space
$\pol(E,\lambda)=\mu_E^{-1}(0)/SO(3)$ is then the symplectic reduction
$$\pol(E,\lambda) = \prod_{e\in E}\calo(\lambda(e))\modmod{0} SO(3)$$
for the diagonal action of $SO(3)$. This determines the symplectic structure
on $\pol(E,\lambda)$.

The codimension 2-embedding 
\begin{equation}\label{eq2}
V_t:=\mu_A^{-1}(\calo(t))\cap\mu_E^{-1}(0) \hookrightarrow
\mu_A^{-1}(\calo(t))\times \mu_{E-A}^{-1}(\calo(t))
\end{equation}
gives rise to a diffeomorphism
\begin{equation}\label{eq3}
\begin{array}{cccccc}
\big[\,V_t/SO(3)\big)\big]/T_A  & \cong & 
\mu_A^{-1}(\calo(t))/SO(3)\times \mu_{E-A}^{-1}(\calo(t))/SO(3)\\[2pt]
\parallel && \parallel \\[2pt]
\pol(E,\lambda)\displaystyle\modmod{t} T_A && 
\pol (\bar A,\lambda^{A,t})\times 
\pol (\overline{E-A},\lambda^{E-A,t}).
\end{array}\end{equation}
As the embedding \eqref{eq2} is the restriction of the obvious symplectomorphism
\begin{equation}\label{eq1}
\prod_{e\in E}\calo(\lambda(e))\ \cong\ \prod_{e\in A}\calo(\lambda(e))\,\times
\!\prod_{e\in E-A}\!\calo(\lambda(e)).
\end{equation}
and as all group actions preserve the symplectic forms, 
the diffeomorphism \eqref{eq3} is a symplectomorphism. \cqfd

\begin{Proposition} \label{dim} Let $\cali$ be a bending 
set for $\pol(E,\lambda)$.  Then
$$\dim T_\cali \leq |E|- \max \{3,|\calm_\cali|\}$$
with equality if and only if $\cali$ is full.
\end{Proposition}

\preu  By Lemma \ref{filling}, it is enough to prove the
formula  when $\cali$ is full. 
We proceed by induction on the number of elements of $\cali$ which are
not singletons. If there are none, then $\dim T_\cali = 0 =
|E|-|E|$ and the formula holds true (recall that $|E|\geq 3$ since we
suppose that $\pol(E,\lambda)\not =\emptyset$). Otherwise, as $\cali$ is full,
there is $A\in\cali$ with $|A|=2$. 

If $|E|=3$, the formula holds true (the $0$-torus, being a quotient of $\bbr^0$, is of
dimension $0$). We may then assume that $|E|\geq 4$.

The map $f_A : \pol(E,\lambda)\to \bbr$ is a moment map for the bending circle
$T_A$. As $|E|\geq 4$, it is not constant. Let $s$ be a regular value of $f_A$
($s>0$ since $A$ is lopsided). By Proposition \ref{redpoly}, the symplectic reduction
of $\pol (E,\lambda)$ at $s$ is a generic polygon space with $|E|-1$ edges.
By Part b) of Remark \ref{redpolyrem}, the bending set 
$\cali$ coinduces a bending set $\bar\cali$ 
for $\bar\lambda$ which is full. The number of non-singletons 
elements of $\bar\cali$ is one less than that of $\cali$. 
By induction hypothesis, one has
$$\dim T_{\bar\cali} =  |E|-1 - \max \{3,|\calm_{\bar\cali}|\} \ .$$
As $\dim T_{\cali} = \dim T_{\bar\cali} +1$ and $\calm_{\bar\cali}=
\calm_{\cali}$, one gets the required expression for $\dim T_\cali$.  \cqfd

\section{Maximal bending tori} \label{mbt}

In this section, we study the poset of bending tori.
Let $\calk$ and $\call$ be two partitions of $E$. We say that $\call$ is {\it
coarser} than $\calk$ if each element of $\call$ is a union 
of elements of $\calk$.

\begin{Theorem} \label{th1rel}  Let $\cali$ be a bending set for
$\pol (E,\lambda)$. Let $N(\lambda ,\cali )$ be the minimal number of
\good\ subsets which are  necessary for
a partition of $E$ which is coarser than $\calm_\cali$. Then,
the maximal dimension $n(\lambda,\cali)$ of a bending torus for $\pol (E,\lambda )$ 
containing $T_\cali$ is
$$n(\lambda ,\cali)= |E|-\max \{3,N(\lambda ,\cali)\}\ .$$
\end{Theorem}

\preu Let $T$ be a bending torus containing $T_\cali$.
By Section \ref{prel}, $T=T_\calj$ for a bending set $\calj$. 
By Lemma \ref{commu}, the
partition $\calm_\calj$ is coarser than $\calm_\cali$.
By \ref{dim}, one has
$$\dim T_\calj \leq |E|-\max \{3,|\calm_\calj|\}\leq
|E|-\max \{3,N(\lambda ,\cali )\}$$
and therefore 
$$n(\lambda ,\cali)\leq |E|-\max \{3,N(\lambda ,\cali)\} .$$

Conversely, let $\calj_0$ be a partition of $E$ into \good\
subsets, coarser than $\calm_\cali$, with $N(\lambda,\cali)$ elements. 
Let $\calj:=\calj_0 \cup \cali$. One
check easily that $\calj$ is a bending set.
Let $\hat\calj$ be a full bending set associated to $\calj$ as in Lemma
\ref{filling}. One has $\calm_{\hat\calj}=\calj_0$ and, by Proposition \ref{dim},
one has,
$$n(\lambda ,\cali) \geq \dim T_{\hat\calj} =
|E|-\max \{3,N(\lambda ,\calj)\}\ . \cqfd $$

As a corollary, we obtain Theorem A of the introduction:

\begin{Theorem}[Theorem A] \label{th1}  Let $N(\lambda )$ 
be the minimal number of
\good\ subsets which are  necessary for a partition of $E$.  
Then the maximal dimension of a bending torus 
for $\pol (E,\lambda)$ is $|E|-\max \{3,N(\lambda )\}$.
\end{Theorem}

\preu  Set $\cali$ be the sets of singletons of $E$ in the statement of Theorem
\ref{th1rel}.  \cqfd

We now give a characterization of the maximal 
bending tori which will be used later.
We can restrict our attention to those  
$T_\cali$, for $\cali$ a full bending set, 
whose dimension is less than $|E|-3$
(the maximal possible dimension of a 
Hamiltonian torus of $\pol (E,\lambda)$).
  
\begin{Proposition}\label{condequi}  Let $\cali$ 
be a full bending set so that
$\dim T_\cali < |E|-3$. 
Then, $T_\cali$ is a maximal bending torus iff
$$\displaystyle\bigcap_{J\in\calm_\calj} 
{\rm Image}(f_J) \not = \emptyset$$
\end{Proposition}

\preu  Observe that $T_\cali$ is a maximal bending torus if and only if for
each pair $I,J\in\calm_\cali$, one has  
${\rm Image}(f_I)\cap {\rm Image}(f_J) \not = \emptyset$
($I\cup J$ is not lopsided). The condition of 
Proposition \ref{condequi} is {\it a priori} 
stronger than that but in fact equivalent, 
thanks to the following lemma.

\begin{Lemma} Let $A_0,\dots ,A_n$ be intervals of the real line. If
$A_i\cap A_j \not = \emptyset$ for all $i,j$, then
$A_1\cap\cdots \cap A_n \not = \emptyset$.
\end{Lemma}

\preu  By induction on $n$, starting with $n=2$. The condition
$A_i\cap A_j \not = \emptyset$ for all $i,j$ implies that
$A:=A_1\cup \cdots \cup A_n$ is connected and hence is an interval.
The set $\cala:=\{A_0,\dots ,A_n \}$ is an acyclic covering of $A$
and therefore its nerve $\caln (\cala )$ can be used to compute
the cohomology of $A$: $H^*(A)=H^*(\caln (\cala ))$. By induction
hypothesis, the simplicial set $\caln (\cala )$ contains the $n-1$ skeleton
of the simplex $\Delta^n$. As $H^{n-1}(A)=0$, $\,\caln (\cala )$
must contain $\Delta^n$ which is to say
$A_1\cap\cdots\cap A_n \not = \emptyset$.  \cqfd

\section{Maximal Hamiltonian tori} \label{mht}

We start with an important special case 
which illustrate the technique: the almost 
regular pentagon. A function
$\lambda : \{1,\dots ,5\}\to \bbr_+$
is called the length function of an {\it almost regular pentagon} if
$\lambda(i)=1$ for $i=1,\dots ,4$ and $1<\lambda(5) <2$. In this case,
$\dim \pol (E,\lambda)=4$.

\begin{Proposition} \label{pentagon} Let  
$\lambda : \{1,\dots ,5\}\to \bbr_+$ be a length
function of an almost regular pentagon. Then, the maximal bending tori of 
$\pol (E,\lambda)$, which are 1-dimensional, are maximal Hamiltonian tori.
\end{Proposition}

\preu  The maximal lopsided subset of $E$ are of the form $\{k,5\}$. Therefore,
all maximal bending tori are of dimension 1. 
Since they are all of the same form, it is enough
to prove Proposition \ref{pentagon} for one of them, say $T_\cali$
with  $\cali:=\{ \{1\}, \{2\},\{3\},\{4,5\}\}$. This gives a
Hamiltonian circle action with moment map 
$f:=f_{\{4,5\}}= |\rho(4)+\rho(5)|$. By Lemma \ref{critfi}, 
this map has three critical values:
\decale{a)} The two extremals $z=\lambda(5)-1$ and 
$z=\lambda(5)+1$ are of course critical values. In both
cases, the critical set is a 2-sphere, the configuration spaces of the
quadrilateral with side length $(1,1,1,z)$.
\decale{b)} the value $1$ for which the critical 
set consists of three points, namely
the configurations $\rho : \{1,\dots ,5\}\to \bbr^3$ given by one of 
the line of equations below
$$\begin{array}{rcccccccc}
-\rho(1) & = & \rho(2) & = & \rho(3) = -\rho(4) - \rho(5), \\
\rho(1) & = & -\rho(2) & = & \rho(3) -\rho(4) - \rho(5) \hbox{ or }  \\
\rho(1) & = & \rho(2) & = & -\rho(3) -\rho(4) - \rho(5).  \\ 
\end{array}$$

The proof then follows from the lemma below.

\begin{Lemma} \label{3pts} 
Let $\mu:M\to\bbr^{m-1}$ be the moment map for a Hamiltonian action of
of $T^{{m-1}}$ on a compact symplectic manifold $M^{2m}$. 
Denote by ${\rm Crit\,}\mu \subset M$ the set of critical points of
$\mu$. Suppose that there is 
a point $\delta$ in the interior of 
the moment polytope $\mu (M)$ such
that $\mu^{-1}(\delta) \cap {\rm Crit\,}\mu$ 
has at least 3 connected
components. Then the action does not extend 
to an effective Hamiltonian action of a $m$-torus.
\end{Lemma}

\preu Suppose that $T$ extends to a Hamiltonian action of $T\times S^1$
with moment map $\Phi : \pol (E,\lambda) \to \bbr^n$. 
Then the moment map $f$ 
is the composition of $\Phi$ with the projection $\bbr^n\to\bbr$ onto
the last coordinate. Additionally, this action, 
being effective, would make $Pol(\lambda)$ a
symplectic toric manifold. Thus, $\Phi(\rho)$ 
are distinct points on the boundary 
of the moment polytope 
$\phi(\pol (E,\lambda))$ (see \cite{De}), which all project to $1$. 
As at most two points of this boundary can project onto 
one point of $\bbr$, we get a contradiction.  \cqfd

The rest of this section is devoted to 
the proof of our second main result:

\begin{Theorem}[Theorem B] \label{mht1}  
Let $T$ be a bending torus 
of $\pol (E,\lambda)$ of
dimension $\geq |E|-5$. Then $T$ is a maximal 
Hamiltonian torus if and only if it is a
maximal bending torus.
\end{Theorem}

We only need to 
prove Theorem B in the cases ${\rm dim\,}T=|E|-4$ and $|E|-5$,
since it is obvious for ${\rm dim\,}T=|E|-3$.

\paragraph{Proof for ${\rm dim\,}T=|E|-4$ :}
Let $\cali$ be a bending set so that 
$T_\cali $ is a maximal bending torus
of dimension $|E|-4$.
We suppose that there is a Hamiltonian circle $S^1$
commuting with $T_\cali$; we shall prove that the resulting action of 
$\widehat T:=T_\cali\times S^1$ is not effective.

Let $f_\cali :\pol (E,\lambda)\to\bbr^{\cali}$ be
the product map $f_\cali:=\prod_{A\in\cali}f_A$. This is
a moment map for the action of $T_\cali$. Its image 
$\Delta$ is a convex polytope of dimension
$|E|-4$. Let $\mu$ be the composition of $f_\cali$ with the projection
to the affine space spaned by $\Delta$ (the ``essential" moment map).

By Proposition \ref{dim}, $\cali$ is full and has 4 maximal elements:
$\calm_\cali = \{I,J,K,L\}$. 
By Proposition \ref{condequi}, there exists a point $c$ 
in the intersection of the images
of $f_I$, $f_J$, $f_K$ and $f_L$. 
The proof divides into 3 cases :

\vskip .2 truecm\noindent\it Case a) : \ \rm
Suppose that $c$ is in the interior of each image. Then
$\vec c:=(c,c,c,c)$ belongs to the interior of the image of the
product map $f:=f_I\times f_J\times f_K\times f_L :\pol (E,\lambda)\to\bbr^4$.
This product map is the composition of $\mu$
with the projection to $\bbr^{\calm_\cali}$.
Hence, there exists $\delta$ in the interior of $\Delta$ which projects
to $\vec c$.

For any  $\rho\in \widetilde\pol(E,\lambda)$ such that $\mu(\rho)=\delta$, there exists
$R_I,R_J,R_K,R_L\in SO(3)$ such that  
$$R_I(\rho_I) = R_J(\rho_J) = -R_K(\rho_K) = -R_L(\rho_L).$$
Then the configuration $\rho'$ defined by 
$$\rho'(e):=R_I(\rho(e)) \hbox{ if } e\in I\ ,\ \rho'(e):=R_J(\rho(e)) \hbox{ if } e\in J, 
\hbox{ etc. }$$
also satisfies $\mu(\rho')=\delta$ and moreover  
$\rho'_I=\rho'_J=-\rho'_K=-\rho'_L$.
This implies that $\rho'$ is a critical point for the function
$h:=f_I+ f_J-f_K- f_L$ and hence for $\mu$. Indeed, the 
Hamiltonian flow of $h$ would be a global rotation around the 
axis $\rho_I$, and therefore induces the identity on $\pol (E,\lambda)$.

Similarly, one constructs critical configurations in $\mu^{-1}(\delta)$ with 
$\rho_I=-\rho_J=\rho_K=-\rho_L$ and $\rho_I=-\rho_J=-\rho_K=\rho_L$.
By lemma \ref{3pts}, this completes the first case.

\vskip .2 truecm\noindent\it Case b) : \ \rm
the  argument of Case a) works as well if $c$ 
is in the interior of the image $f_A$ for each
$A\in \calm_\cali$ which is not a singleton (by genericity of $\lambda$, 
there exists at least one such element). 

\vskip .2 truecm\noindent\it Case c) : \ \rm
in the general case, there may be some set $A\in\calm_\cali$, 
such that $c$ is in the boundary of the image of 
$f_A$. Let $\calm'\subset \calm_\cali$ be the set of such $A$'s and
let $\bar \calm'$ be the partition of $E$ generated by
$\calm'$ (formed by the elements of $\calm'$ and the singletons).
Call $\cali'$ the largest sub-poset of $\cali$ 
so that $\calm_{\cali'}=\bar\calm'$; this is a full bending set.

In this case, $\bar P:=f^{-1}(\vec c\, )$ is a symplectic submanifold of $\pol (E,\lambda)$
on which $T_{\cali'}$ acts trivially. As $\bar P$ coincides with the result of 
successive symplectic reductions at $c$ for the various $f_A$ with $A\in\calm'$,
it is, by Proposition \ref{redpoly},
symplectomorphic to the polygon space 
$\pol(\bar\calm',\bar\lambda)$, where 
$$\bar\lambda (\{e\})=\lambda(e) \ \hbox{ and } \ 
\bar\lambda (A) = c \  \hbox{ if } A\in\calm'$$
The bending torus $T_\cali$ acts on $\bar P$, giving rise to a bending torus
$T_{\bar I}$ isomorphic to $T_\cali/T_{\cali'}$. Observe that
$\bar I$ has 4 maximal elements and that we are in Case b).  Therefore,
$T_{\bar I}$ is a maximal Hamiltonian torus and the induced action of 
$\widehat T$ on $\bar P$ has a kernel of dimension strictly larger than that of $T_{\cali'}$.
Therefore, as
$$  {\rm dim\,}\pol (E,\lambda) - {\rm dim\,}\bar P = 
2\big( \sum_{A\in\calm'} |A| - |\calm'|\big) = 2\, {\rm dim\,} T_{\cali'},$$
there is a circle in $\widehat T$ acting trivially on a tubular neighborhood
of $\bar P$. Hence, by the generic orbit type theorem \cite[\S\ 2.2]{Au},
the action of $\widehat T$ on $\pol(E,\lambda)$ is not effective.  \cqfd

\paragraph{Proof for ${\rm dim\,}T=|E|-5$ : }
Let $\cali$ be a bending set so that 
$T_\cali $ is a maximal bending torus
of dimension $|E|-5$. 
We suppose that there is a Hamiltonian circle $S^1$
commuting with $T_\cali$ and we shall prove that the resulting action of 
$\widehat T:=T_\cali\times S^1$ is not effective.

Let $\mu : \pol(E,\lambda)\to\bbr^{|E|-5}$ be
the essential moment map, defined 
as in the proof for ${\rm dim\,}T = |E|-4$, and let and $\Delta$ be the image of $\mu$. 
Let $\widehat\mu :\pol(E,\lambda)\to\Delta\times\bbr$ be a moment map for the action of 
$\widehat T$ with first component equal to $\mu$ and let $\widehat\Delta$
be the image of $\widehat\mu$.

By Proposition \ref{dim}, $\calm_\cali$ has 5 elements. By Proposition
\ref{condequi}, there exists a point $c$ in the intersection of the images
of $f_A$ for $A\in\calm_\cali$. The proof divides into several cases :

\vskip .2 truecm\noindent\it Case 1) : \ \rm Suppose that $|E|=5$. Then $T_I$ is
of dimension 0 and we have to know that a maximal Hamiltonian torus
for a regular pentagon space is also of dimension 0. This is the contents of
\cite[Theorem 3.2]{HK2}.

\vskip .2 truecm\noindent\it Case 2) : \ \rm
Suppose that each $A\in\calm_\cali$ contains exactly 2 elements 
(hence $|E|=10$) and $c$ is in 
the interior of the image of $f_A$.
This implies that
$\vec c:=(c,c,c,c,c)$ is a regular value of $\mu$. 
The reduction $Q$ of $\pol (E,\lambda)$ at $\vec c$ is then symplectomorphic
to a regular pentagon space (apply Proposition \ref{redpoly} five times). 
The induced Hamiltonian action of $\widehat T$ on $Q$ is then trivial by Case 1). 
This implies that
the image of the differential $D\widehat\mu$ 
at any point of $\mu^{-1}(\vec c\, )$ is parallel to $\Delta\times\{0\}$. 
By convexity, we deduce that $\widehat\Delta$ and $\Delta$ have the same dimension 
and therefore the action of $\widehat T$ is not effective.

\vskip .2 truecm\noindent\it Case 3) : \ \rm
The argument of Case 2) works as well if each 
$A\in\calm_\cali$ has $\leq 2$
elements and $c$ is in 
the interior of the image of $f_A$ when $|A|=2$. Also, if there are 
sets $A\in\calm_\cali$ with $|A|=2$ and $c$ is in the boundary of 
the image of $f_A$, one proceeds as in Case c) of the proof for 
${\rm dim\,}T_\cali = |E|-4$ to deduce that the action of $\widehat T$ is not
effective. Thus, we are able to prove our result when all the elements of $\calm_\cali$
are either singletons or doubletons.

\vskip .2 truecm\noindent\it General case) : \ \rm For $A\in\calm_\cali$, let
$k_A:=\max \{0,|A|-2\}$ and $k:=\sum_{A\in\calm_\cali} k_A$. The proof goes by induction
on $k$, the case $k=0$ being established in Case 3). If $k>0$, let $A\in\calm_\cali$
such that $|A|\geq 3$. If $c$ lies in the boundary of the image of $f_A$, one 
proceeds as in Case c) of the proof for 
${\rm dim\,}T_\cali = |E|-4$ to deduce that the action of $\widehat T$ is not
effective (using the induction hypothesis). Otherwise, as $\cali$ is full, there exists $B\in\cali$ such that
$|B|=2$, $B\subset A$ and $f_B(f_A^{-1}(c))$ is an interval of positive length.
It contains an open interval $J$ of regular values of $f_B$. For $t\in J$, the reduction
of $\pol(E,\lambda)$ for the action of the Hamiltonian circle with moment map $f_B$
is, by Proposition \ref{redpoly}, symplectomorphic to an $(|E|-1)$-gon space $\bar P$.
The bending torus $T_\cali$ descends to a bending torus $T_{\bar\cali}$ for $\bar P$.
One has $\calm_{\bar\cali}=\calm_\cali$ and $\bar k = k-1$. By induction hypothesis,
$T_{\bar\cali}$ is a maximal Hamiltonian torus. This implies that each point
of $f_B^{-1}(t)$ has a stabilizer of positive dimension for the action of $\widehat T$.
This holds true for all $t\in J$, therefore for an open set of $\pol(E,\lambda)$. By the
generic orbit type theorem \cite[\S\ 2.2]{Au}, this implies that the action of $\widehat T$
on $\pol(E,\lambda)$ is not effective.  \cqfd

\section{Examples}\label{ex}

\noindent{\sc Notations : } When $E=\{1,\dots ,n\}$, we describe $\pol (E,\lambda)$
by writing the values of $\lambda$. For instance, $\pol(1,1,1,2)$ stands for
$\pol(\{1,2,3,4\},\lambda)$ with $\lambda(1)=\lambda(2)=\lambda(3)=1$ and
$\lambda(4)=2$. A bending set is described by listing 
its elements which are not singletons and labeling the edges by their length.

\begin{ccote}\label{twolong}
{The ``two long edge" case : \ }\rm Suppose that the set of edges
$E$ contains two elements $a,b$ such that
$$\lambda(a)+\lambda(b) > \sum_{e\in E-\{a,b\}}\lambda(e) \ .$$
Then $E$ is the disjoint union of $E_a$ and $E_b$ so that $E_a$ is lopsided
with longest edge $a$ and
$E_b$ is lopsided with longest edge $b$. One then has
$N(\lambda)=2$ and, by Theorem \ref{th1rel}, $\pol (E,\lambda)$ admits a bending
torus of dimension $|E|-3$. In particular, $\pol (E,\lambda)$ is a toric
manifold.
\end{ccote}

\begin{ccote}
{Almost regular pentagon :\ }\rm
The almost regular pentagon $\pol (1,1,1,1,a)$ with $1<a<2$
(or $0<a<1$)
is a very important special case, already used in Proposition \ref{pentagon}.
Notice  $\pol (E,\lambda)$ is diffeomorphic to
$\bbc P^2 \, \sharp \, 4\, \overline{\bbc P^2}$
(see \cite[Example 10.4]{HK1}). 

We used the result of [HK2] that the regular pentagon space admits no non-trivial 
circle action. This is not known for regular polygon spaces with more edges.
Nor it is known whether an almost regular pentagon space
is diffeomorphic to a toric manifold.
\end{ccote}

\begin{ccote} \label{toridiffdim}
{Hamiltonian tori of different dimensions :\ }\rm
Consider a generic pentagon space of the form $P_{a,b}:=\pol(1,1,1,a,b)$
with $a\not = 1 \not = b$ and $0<a-b<1<a+b$. 
The bending circle $\{a,b\}$ is a maximal Hamiltonian torus by 
Proposition \ref{condequi} and \ref{mht1}. However, $\pol (1,1,1,a,b)$ is
a toric manifold by the bending tori $T_\cali$ of the form
$\cali:=\{\{1,a\},\{1,b\}\}$. 
In this example, one sees that maximal bending tori, as well
as maximal Hamiltonian tori, are not all of the same dimension. 

The moment polytope for $T_\cali$ shows that $P_{a,b}$ is
diffeomorphic to $\bbc P^2 \, \sharp \, 4\, \overline{\bbc P^2}$
if $a+b<3$ and to
$\bbc P^2 \, \sharp \, 3\, \overline{\bbc P^2}$
if $a+b>3$ (the case $a+b=3$ is not generic).
It is known that the other pentagon spaces are $4$-manifolds
with second Betti number $<3$. For them, any Hamiltonian circle action
extends to a toric action by \cite[Th.\,1]{Ka}.

An example with maximal Hamiltonian tori of 3 different dimensions is
provided by the heptagon spaces
$\pol (1,1,2,2,3,3,3)$ (it is generic since lengths are integral and the 
perimeter is odd). The 3 bending sets with maximal 
(non-singleton) elements of the form
$$\{\{2,1\},\{2,1\}\} \quad ,\quad 
\{\{2,1\},\{3,1\},\{3,2\}\} \quad ,\quad 
\{\{3,1,1\},\{3,2\},\{3,2\}\}$$
determine maximal Hamiltonian tori of dimension respectively $2,3$ and $4$.
Observe that the bending circle $\{3,2\}$ is contained in two maximal tori
of different dimension. 

Examples in higher dimension can be constructed by adding ``little edges"
to the previous one, for instance
the ($7+m$)-gon space 
$$\pol(1,1,2,2,3,3,3,1/2,1/4,\dots ,1/2^m).$$
It admit full bending sets with maximal (non-singleton) elements of the form \goodbreak
\begin{itemize}
\item $\{\{2,1\},\{2,1\},\{3,1/2,1/4,\dots ,1/2^m\}\}$ 
\item $\{\{2,1\},\{3,1\},\{3,2\},\{3,1/2,1/4,\dots ,1/2^m\}\}$
\item $\{\{3,1,1\},\{3,2\},\{3,2\},\{3,1/2,1/4,\dots ,1/2^m\}\}$
\end{itemize}
which determine maximal Hamiltonian tori of dimension respectively 
$m+2,m+3$ and $m+4$.
\end{ccote}

\begin{ccote}\label{notbend}\rm
Let $T_1$ and $T_2$ be two Hamiltonian tori of dimension $n$
for a symplectic manifold $M^{2n}$.
Choose isomorphisms
${\rm Lie}(T_1)^*\approx\bbr^n\approx {\rm Lie}(T_2)^*$.
the moment polytopes $\Delta_1$ and $\Delta_2$ of the two actions
are in $\bbr^n$. By Delzant's theorem, $T_1$ is conjugate to $T_2$
in the group $\cals(M)$ of sympectomorphism of $M$ if and only if
the moment polytopes $\Delta(T_i)$ satisfy
$\Delta(T_2)=\psi(\Delta(T_1))$ where $\psi$ is a
composition of translations and transformations in $GL(\bbz^n)$.

Consider the pentagon space $P:=\pol(1,a,c,c,c)$,
with $c>a+1>2$.
The two bending tori
$T_1=\{\{c,1\},\{c,a\}\}$ and $T_2=\{\{c,1\},\{c,a,1\}\}$
have moment polytopes

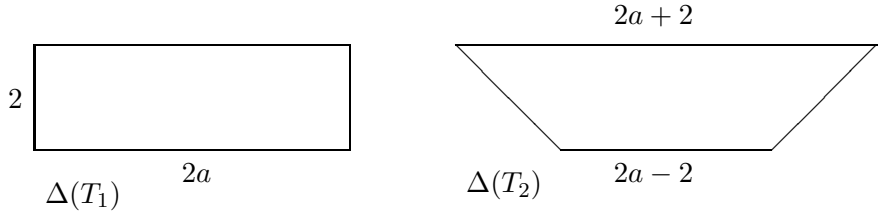
\setlength{\unitlength}{.7mm}
\begin{picture}(10,20)(25,15)
\put(20,20){\line(1,0){60}}
\put(20,0){\line(1,0){60}}
\put(20,0){\line(0,1){20}}
\put(80,0){\line(0,1){20}}
\put(48,-6){$2a$}
\put(15,8){$2$}
\put(22,-10){$\Delta(T_1)$}

\put(100,20){\line(1,0){80}}
\put(120,0){\line(1,0){40}}
\put(100,20){\line(1,-1){20}}
\put(160,0){\line(1,1){20}}
\put(130,-6){$2a-2$}
\put(130,24){$2a+2$}
\put(102,-8){$\Delta(T_2)$}
\end{picture}\vskip 25mm
Therefore, $T_1$ and $T_2$ are not conjugate
in in the group $\cals (P)$.
One can check that any other bending torus is
conjugate to either $T_1$ or $T_2$.

On the other hand, the polytope $\Delta(T_1)$
shows that $P$ is symplectomorphic to
$(S^2\times S^2,\omega_1+a\omega_2)$, where
$\omega_1$ and $\omega_2$ are the pull back of
the standard area form on $S^2$ via the two projection maps.
By \cite[Th.\,2]{Ka}, the number of conjugacy classes of maximal
Hamiltonian tori is equal to $[a]$, the smallest integer
greater than or equal to $a$. This proves the following

\begin{Proposition}\label{pnotbend}
If $c>a+1>3$, then $\pol(1,a,c,c,c)$ admits
Hamiltonian tori which are not conjugate to a bending torus.
\end{Proposition}
\end{ccote}

\begin{ccote}\label{Lerman}  \rm
Let $(M,\omega)$ be a simply connected symplectic manifold
such that $[\omega]\in H^2(M;\bbr)$ is integral.
Then there exists
a principal circle bundle $S^1\to Q\to M$
with Euler class $[\omega]$ and $Q$ carries a natural contact
distribution by a theorem of Boothby and Wang
\cite[Th.\,3]{BW}. In \cite[Th.\,1]{Le}, E. Lerman
recently proved that maximal Hamiltonian tori
in $M$ (of dimension $k$) give rise to
maximal tori (of dimension $k+1$) in the group of
diffeomorphism of $Q$ preserving the contact distribution.

By \cite[Prop.\,6.5]{HK1}, the symplectic form on $\pol(E,\lambda)$
is integral when, for example, $\lambda$ takes integral values.
Then, our examples in \ref{toridiffdim} give rise to
contact manifolds with maximal tori of
different dimensions in their group of
contactomorphisms (see \cite[Example 2]{Le}).
\end{ccote}

\small
\noindent Jean-Claude HAUSMANN, Math\'ematiques-Universit\'e B.P. 240 \\ CH-1211 Gen\`eve
24, Suisse.\\ hausmann@math.unige.ch
\vskip 0.4 truecm
\noindent Susan TOLMAN, Department of Mathematics,
University at Illinois at Urbana-Champaign
Urbana, IL 61801, USA\\ stolman@math.uiuc.edu

\begin{thebibliography}{HK2}
\small

\bibitem[Au]{Au} Audin M.
The topology of torus actions on symplectic manifolds.
{\em Birkh\"auser} (1991).
\bibitem[BW]{BW} W.M. Boothby and H.C. Wang, On contact manifolds,
{\em Ann.\ of Math.\ (2)} {\bf 68} (1958), 721--734.
\bibitem[De]{De} Delzant T.
Hamiltoniens p\'eriodiques et image convexe de l'application moment.
{\em Bull. Soc. Math. France}
{\bf 116} (1988), {315--339}
\bibitem[HK1]{HK1}
Hausmann, J-C. \& Knutson A. The cohomology ring of polygon spaces.
Grasmannians. {\em Annales de l'Institut Fourier} (1998) 281-321.
\bibitem[HK2]{HK2} Hausmann, J-C. \& Knutson A.
A limit of toric symplectic forms that has no 
periodic Hamiltonians.
{\em GAFA, Geom. funct. anal.} 10 (2000) 556--562.
\bibitem[Ha]{Ha}
Hausmann, J-C.
Sur la topologie des bras articul\'es.
{\em In ``Algebraic Topology, Poznan", Springer Lectures Notes}
{\bf 1474} (1989), 146--159.
\bibitem[KM]{KM}
Kapovich, M. \& Millson, J. {The symplectic geometry of polygons in Euclidean
space.} {\em J. of Diff. Geometry} {\bf 44} (1996), 479--513.
\bibitem[Ka]{Ka}
Karshon, Y.
Maximal tori in the symplectomorphism groups
of Hirzebruch surfaces.
Preprint, http://www.ma.huji.ac.il/~karshon/papers
\bibitem[Kl]{Kl}
Klyachko, A. Spatial polygons and stable configurations
of points in the projective line.
{\em in: Algebraic geometry and its applications 
(Yaroslavl, 1992), Aspects
Math., Vieweg, Braunschweig} (1994) 67--84.
\bibitem[Le]{Le}
Lerman, E.
On maximal tori in the contactomorphism groups of regular contact
manifolds. Preprint, 2002.
\end{thebibliography}
\end{document}